\documentclass{article}
\usepackage{times}
\usepackage{amsmath}
\usepackage{amsfonts} 
\usepackage{lipsum}
\usepackage{ntheorem}
\theoremstyle{plain}
\newtheorem*{theorem 1}{The Schur 1916 Theorem 3}

\begin{document}

\title{Off-Diagonal Continuous Rado Numbers for $x_1 + x_2 + \cdots + x_k = x_0$}
\author{Dr. Donald Vestal and Jonathan Sax,\\
 Department of Mathematics \& Statistics,\\
  South Dakota State University,\\
  Brookings,\\
  South Dakota,\\
  \texttt{Donald.Vestal@sdstate.edu} \\
 \texttt{Jonathan.Sax@jacks.sdstate.edu}}
\date{\today}
\maketitle

\begin{abstract}
\par\noindent
In 2001, Robertson and Schaal found the 2-color off-diagonal generalized Schur numbers: for two positive integers $k$ and $l$, they determined the smallest positive integer $S = S(k, l)$ such that for any coloring of the integers from 1 to $S$ using red and blue, there must be a red solution to the equation $x_1 + x_2 + \dots + x_k = x_0$ or a blue solution to the equation $x_1 + x_2 + \dots + x_l = x_0$. We extend this result to find the continuous version: for two positive integers $k$ and $l$, we find the smallest real number $S = S_\mathbb{R} (k, l)$ such that for any coloring of the real numbers from 1 to $S$ using red and blue, there must be a red solution to the equation $x_1 + x_2 + \dots + x_k = x_0$ or a blue solution to the equation $x_1 + x_2 + \dots + x_l = x_0$. 

\end{abstract}

\section{Introduction}

\par
Ramsey Theory is a branch of mathematics that finds the conditions under which order must appear. Generally, the idea is that when a structure is ``large'' enough, we can guarantee that a certain substructure appears. There are many Ramsey theory results in graph theory, but in this article, we will focus on Ramsey theory applied to the natural numbers, and then extend this to the positive real numbers.

\par In 1916 (predating Ramsey's 1928 paper that is often considered the beginning of Ramsey Theory \cite{Ramsey}), Isaai Schur proved the following.

\par\noindent
\textbf{Schur's Theorem}  For any positive integer $t$ there is an integer $S=S(t)$ such that any $t$-coloring of the set $\left\{ 1, 2, \dots , S \right\}$ must result in three integers $x, y, z$, all of which are colored the same and such that $x + y = z$.

\par 
Following this, Richard Rado (a student of Isaai Schur) furthered the work of his mentor to find a similar results \cite{Rado} for linear equations and, more generally, systems of linear equations. So, under certain conditions, for a given system of linear equations, if the positive integers are colored using $t$ colors, then we can guarantee a monochromatic solution to this system; in other words, a solution in which all of the positive integers in that solution have the same color. 

\par
Rado's result concerned the existence of a monochromatic solution; however, for a given equation (or system of equations) $\mathcal{E}$ and a given number $t$ of colors, if there must eventually be a monochromatic solution, then there must be a smallest positive integer by which we can guarantee this monochromatic solution. This smallest number is often referred to as ``the $t$-color Rado number for $\mathcal{E}$'', except in the special case of the equation $x + y = z$, where we simply refer to it as the $t$-color Schur number, denoted $S(t)$, due to Schur's theorem above. The known Schur numbers are: $S(1) = 2, S(2) = 5, S(3) = 14, S(4) = 45,$ and $S(5) = 161$ \cite{OEIS}. 

\par
In 1982, Beutelspacher and Brestovansky \cite{BB} generalized Schur's result by proving the following.

\par\noindent
\textbf{Beutelspacher and Brestovansky} For any positive integer $k$, the 2-color Rado number for the equation $x_1 + x_2 + \dots + x_k = x_0$ is $k^2 + k - 1$.

\par
There have been many families of equations that have been studied; below we give a list of some examples, with references.

\begin{tabbing}
\par\noindent \qquad \= \cite{SST} \quad \= $x_1 + x_2 + \cdots + x_{m-1} + c = x_m$ where $c \geq 0$ \\
\> \cite{KS} \> $x_1 + x_2 + \cdots + x_{m-1} + c = x_m$ where $c < 0$ \\
\> \cite{SV} \> $x_1 + x_2 + \cdots + x_{m-1} = 2 x_m$ \\
\> \cite{GS} \> $a_1 x_1 + a_2 x_2 + \cdots + a_m x_m = x_0$ \\
\> \cite{DSV} \> $x + y^n = z$ 
\end{tabbing}

\par
For this article, there are two results that we will focus on. The first is a continuous version of the Rado number, in which we color the real numbers in an interval (not just the integers). In \cite{BH}, Brady and Haas prove the following: 
\par\noindent
\textbf{Brady and Haas} For a given positive real number $\gamma$ and a positive integer $m$, suppose that $c$ is a real number which satisfies $-c < \gamma \left( m-1 \right)$.  Let $R = R \left( \gamma , m , c \right)$ denote the least real number such that for any 2-coloring of the real numbers in the interval $\left[ \gamma , R \right]$, there must be a monochromatic solution to the equation $x_1 + x_2 + \dots + x_m + c = x_0$. Then $R \left( \gamma , m , c \right) = \gamma \left( m^2 + m - 1 \right) + \left( m + 2 \right) c$.

\par
The second is an ``off-diagonal'' result from Robertson and Schaal \cite{RS}, in which they focus on the following generalization of Beutelspacher and Brestovansky: 
\par\noindent
\textbf{Robertson and Schaal} Given two positive integers $k$ and $l$ with $2 \leq k \leq l$, let $S = S(k,l)$ denote the smallest integer such that for any 2-coloring of the integers in $\left\{ 1, 2, \dots S \right\}$ using red and blue, there must be a red solution to the equation $x_1 + x_2 + \cdots + x_k = x_0$ or a blue solution to the equation $x_1 + x_2 + \cdots + x_l = x_0$.  Then
\par\noindent
\[
    S(k, l) =
\begin{cases}
    3l-1,& \text{if  } k = 2 \text{ and } l \geq 2 \text{ and is even,}\\
    3l-2,& \text{if  } k = 2 \text{ and } l \geq 3 \text{ and is odd,}\\
    kl + k - 1,&\text{if } 3\leq k \leq l \text{.}
\end{cases}
\]

\par\noindent 
In the case where $k = l$, this of course gives the result in \cite{BB}, and more specifically for $k=l=2$, this gives the Schur number $S(2)$. 

\par
We will prove a continuous version of this off-diagonal result of Robertson and Schaal: for two positive integers $k$ and $l$ with $2 \leq k \leq l$, we let $S = S_{\mathbb{R}} (k,l)$ denote the smallest real number for which any 2-coloring of the interval $\left[ 1, S \right]$ using red and blue must result in a red solution to the equation $x_1 + x_2 + \cdots + x_k = x_0$ or a blue solution to the equation $x_1 + x_2 + \cdots + x_l = x_0$. 
\par\noindent
\textbf{Theorem 1} If $k$ and $l$ are positive integers with $2 \leq k \leq l$ , then $S_\mathbb{R} (k, l) = kl + k - 1$.

\par
Note that when $3 \leq k \leq l$, the formula for the continuous version matches the discrete formula from \cite{RS}. For the most part, we will only have to establish the coninuous result for $k = 2$, and then use the result from \cite{RS} for the $k \geq 3$ case, except that we will correct one error: in Case III for the upper bound in \cite{RS}, they make the claim that we can assume that 1 is red, since a proof when 1 is blue can be obtained by interchanging red and blue along with $k$ and $l$. This turns out not to be true, but we will adjust the proof to fix that error. The error comes about for the following reason: we are assuming that $k < l$ and that the red solutions involve the equation with fewer variables. When we assume 1 is blue, then the resulting coloring is affected. Specifically, we get a longer coloring if we start by coloring 1 with the color associated with the equation with fewer variables, whereas starting by coloring 1 with the color associated with the equation with more variables results in a shorter coloring.   

\section{The Lower Bound}
For the remainder of this article, let $S = kl + k - 1$. To prove our result, we need to establish the lower bound: there is a coloring of the interval $\left[ 1, S \right)$ which avoids both a red solution to $x_1 + x_2 + \cdots + x_k = x_0$ and a blue solution to $x_1 + x_2 + \cdots + x_l = x_0$. For the remainder of this article, we will use the term ``monochromatic solution" to refer to either a red solution to $x_1 + x_2 + \cdots + x_k = x_0$ and a blue solution to $x_1 + x_2 + \cdots + x_l = x_0$. Thus, for example, a blue solution to the equation $x_1 + x_2 + \cdots + x_k = x_0$ will not be considered a monochromatic solution. 
\par\noindent
\textbf{Lemma 1} The following coloring of the interval $\left[ 1, S \right)$  avoids a red solution to $x_1 + x_2 + \cdots + x_k = x_0$ and a blue solution to $x_1 + x_2 + \cdots + x_l = x_0$:
\par\noindent
\qquad Red: $\left[ 1, k \right) \cup \left[ kl, S \right)$
\par\noindent
\qquad Blue: $\left[ k, kl \right) $.
\par\noindent
Proof. First, there is clearly no blue solution to the equation $x_1 + x_2 + \cdots + x_l = x_0$ since the smallest value the left hand side can take occurs when all variables on the left are equal to $k$, which means on the right, we will have $x_0 \geq kl$, and there are no blue numbers that satisfy this inequality. For the red equation, note that if each variable on the left side satisfies $1 \leq x_i < k$, then adding these up results in $k \leq x_0 < k^2 $; but $k \leq l$, which means $k \leq x_0 < kl$, making $x_0$ colored blue. So at least one variable on the left side of the red equation would have to take on a value in the interval $\left[ kl, S \right)$. But if there is a variable from the left side that lies in that interval, then the smallest value $x_0$ could take would be if the remaining variables on the left were equal to 1, giving $x_0 \geq \left( k - 1 \right) + kl = S$ and thus $x_0$ is not in the interval $\left[ 1, S \right)$.  \hfill{$\square$}

\section{The Upper Bound}

\par
To complete our proof, we need to establish the upper bound: for any coloring of the real numbers in the closed interval $\left[ 1, S \right]$ using red and blue, there must be a red solution to $x_1 + x_2 + \cdots + x_k = x_0$ or a blue solution to $x_1 + x_2 + \cdots + x_l = x_0$. 

\par
We begin with the case $k = 2$. Note that the discrete case in \cite{RS} proves that there will be a solution involving integers between 1 and either $3l-2$ or $3l-1$ (depending on the parity of $l$). But our continuous Rado number $S = kl + k - 1 = 2l + 1$ is smaller since we can use solutions that are not integers. (And thus our proof at some point must involve some non-integer values.)

\par\noindent
\textbf{Lemma 2} If the real numbers in the interval $\left[ 1, 2l+1 \right)$ are colored red and blue, then there must be a red solution to $x_1 + x_2 = x_0$ or a blue solution to $x_1 + x_2 + \cdots + x_l = x_0$.

\par\noindent
Proof. We consider the two possible colorings for the number 1. Assume that we are trying to avoid a red solution to $x_1 + x_2 = x_0$ and a blue solution to $x_1 + x_2 + \cdots + x_l = x_0$.
\par\noindent 
Case 1: 1 is colored red
\par
If 1 is red, then the solution $\left( 1,  1, 2 \right)$ would force 2 to be blue. With the solution $\left( 2, 2, \dots , 2, 2l \right)$, $2l$ would be colored red. With 1 and $l$ colored red, the two solutions $\left( 1, 2l, 2l+1 \right)$ and $\left( 1, 2l-1, 2l \right)$ would make both $2l-1$ and $2l+1$ blue. The two solutions $\left( 2, 2, \dots , 2, \frac{3}{2} , \frac{3}{2} , 2l-1 \right) $ and $\left( 2, 2, \dots , 2, \frac{5}{2} , \frac{5}{2} , 2l+1 \right) $  would then make both $\frac{3}{2}$ and $\frac{5}{2}$ red. But this then results in the red solution $\left( 1, \frac{3}{2} , \frac{5}{2} \right) $.

\par\noindent
Case 2: 1 is colored blue
\par
If 1 is blue, then the solution $\left( 1, 1, \dots , 1, l \right)$ would force $l$ to be red. With the solution $\left( l, l, 2l \right)$, $2l$ would be colored blue, and from the solution $\left( 2, 2, \dots , 2, 2l \right)$, 2 would have to be colored red. With 2 and $2l$ colored red, the two solutions $\left( 2, 2, 4 \right)$ and $\left( 2, l, l+2 \right)$ would make both 4 and $l+2$ blue. The two solutions $\left( 1, 1, \dots , 1, 3, l+2 \right) $ and $\left( 1, 1, \dots , 1, 4 , l+3 \right) $  would then make both 3 and $l+3$ red. But this then results in the red solution $\left( 3, l, l+3 \right) $.  \hfill{$\square$} 

\par
For $k \geq 3$, we can use the proof from \cite{RS}: we have a lower bound coloring for the interval $ \left[ 1, S \right)$, and the argument in \cite{RS} shows that there is an integer solution in the interval $\left[ 1, S \right]$, so we obviously have at least one real solution. However, there is a piece of the argument in \cite{RS} that needs to be fixed. As mentioned eariler, in Case III for the upper bound in \cite{RS}, they make the claim that we can assume that 1 is red, since a proof when 1 is blue can be obtained by interchanging red and blue along with $k$ and $l$. But if the argument is followed with the values of $k$ and $l$ switched, then a solution involving $kl + l - 1$ is used; this is a problem since $kl + l - 1$ is larger than $kl + k - 1$ which is the number that is claimed as the 2-color off-diagonal Rado number. So, we fix this with two more lemmas.

\par\noindent
\textbf{Lemma 3} Suppose $l > k$. If 1 is colored blue, then to avoid a red solution to $x_1 + x_2 + \cdots +x_k = x_0$ and a blue solution to $x_1 + x_2 + \cdots +x_l = x_0$ we must color $k, l,$ and $l+1$ red and $kl$ blue.

\par\noindent
Proof. If 1 is blue, then the solution $\left( 1, 1, \dots , 1, l \right)$ to the blue equation forces $l$ to be colored red. The solution $\left( l, l, \dots , l, kl \right)$ to the red equation then forces $kl$ to be colored blue. And the solution $\left( k, k, \dots , k, kl \right)$ to the blue equation then forces $k$ to be colored red. Finally, the solution $x_1 = x_2 = \cdots = x_{l-k+1} = 1, x_{l-k+2} = \cdots = x_l = l+1, x_0 = kl$ to the blue equation forces $l + 1$ to be colored red. \hfill{$\square$} 

\par\noindent
\textbf{Lemma 4} Suppose $l = k + r$ where $r > 0$. If 1 is colored blue, then 2 must be colored red to avoid a monochromatic solution.

\par\noindent
Proof. Since we are avoiding a monochromatic solution, from Lemma 3, we currently have $k, l, $ and $l + 1$ colored red and 1 and $kl$ colored blue.  Let $b$ be the smallest nonnegative interger with $b \equiv 1-k \ ( \text{mod } r)$. (So $0 \leq b < r$.) Then let 
$$y = k - 1 + \frac{1-k-b}{r},$$
\par\noindent
which will be an integer satisfying $k-2-\frac{k-1}{r} < y \leq k - 1 - \frac{k-1}{r}$. From the solution $x_1 = \cdots = x_{k-y} = k, x_{k-y+1} = \cdots = x_k = l, x_0 = k^2 + \left( r - 1 \right) \left( k - 1\right) - b$, we must have $k^2 + \left( r - 1 \right) \left( k - 1\right) - b$ colored blue in order to avoid a red solution. 
\par
If $b=0$, then we will have a blue solution to the blue equation: 
$$\left( 1, 1, \dots , 1, k^2 + \left( r - 1 \right) \left( k - 1\right) , kl \right)$$
\par\noindent
thus giving us a monochromatic solution and completing our proof of Theorem 2. On the other hand, if $b \neq 0$, then from the solution 
$$x_1 = \cdots = x_{l-b-1} = 1,$$
$$ x_{l-b} = \cdots = x_{l-1} = 2,$$
$$ x_l = k^2 + \left( r - 1 \right) \left( k - 1\right) - b,$$
\par\noindent
and
$$ x_0 = kl$$ 
\par\noindent
we conclude that 2 must be red to avoid a blue solution to the blue equation. \hfill{$\square$} 

\par
With these lemmas, we can prove the existence of a monochromatic solution.

\par\noindent
\textbf{Theorem 2} Suppose $l > k$. If 1 is colored blue and 2 and $l+1$ are colored red, then we must have a monochromatic solution.

\par\noindent
Proof. Since 2 and $l+1$ are red, the solutions 
$$\left( 2, 2, \dots , 2, 2k \right) \text{ and } \left( 2, 2, \dots 2, l+1, 2k+l-1 \right)$$
\par\noindent
would force $2k$ and $2k+l -1$ to be blue in order to avoid a red solution. However, this then results in the blue solution $\left( 1, 1, \dots , 1, 2k, 2k+l-1 \right)$.  \hfill{$\square$} 

\par
This completes the proof for Theorem 1. 

\par
It is worth noting that, due to the fact that our equations are homogeneous, we can extend the theorem to start at any positive real number. For positive integers $k$ and $l$ with $2 \leq k \leq l$ and a positive number $\gamma$, let $S_{\mathbb{R}} (\gamma, k, l)$ denote the smallest real number $S$ such that for any coloring of the real numbers in the interval $\left[ \gamma , S \right]$ using red and blue, there must be a red solution to the equation $x_1 + x_2 + \cdots + x_k = x_0$ or a blue solution to the equation $x_1 + x_2 + \cdots + x_l = x_0$. 

\par\noindent
\textbf{Theorem 3} If $k$ and $l$ are positive integers with $2 \leq k \leq l$ and $\gamma > 0$, then $S_\mathbb{R} (\gamma , k, l) = \gamma kl + \gamma k - \gamma$.

\section{Conclusion}

One thing worth noting in this continuous result is how moving to the real numbers has the effect of ``smoothing'' the formulas for the smaller cases. In \cite{RS}, there is a general formula that holds when $k \geq 3$, but a slightly different result when $k = 2$; when we look at the continuous version, we get a formula that holds for all integers $k \geq 2$. In fact, we could extend this trivially to the case $k = 1$, since the red equation would just be $x_1 = x_0$. The only way to avoid a red solution to this equation would be to color every real number blue. This gives a blue solution to the equation $x_1 + x_2 + \cdots + x_l = x_0 $ once you color $l$. Thus, the Rado number would then be $l$, which agrees with the formula $kl + k - 1$ when $k = 1$.

\par
For future work, there are other families of equations that could be studied for this off-diagonal generalization, including using three or more colors (although that quickly makes things more difficult). But also the continuous versions could be studied, which seems to be a largely untouched area in Ramsey theory.

\end{document}